\documentclass[aop,MSNbibl,dvips]{arximspdf}
\usepackage{mathbh}

\doi{10.1214/12-AOP817} 
\volume{41}
\issue{6}
\pubyear{2013}
\firstpage{3879}
\lastpage{3909}

\makeatletter
\newcommand{\rrvert}{\vert}
\newcommand{\llvert}{\vert}
\newcommand{\eqref}[1]{(\ref{#1})}
\newtheorem{theorem}{Theorem}[section]
\newproclaim{definition}[theorem]{Definition}
\newtheorem{proposition}[theorem]{Proposition}
\newtheorem{lemma}[theorem]{Lemma}

\newproclaim{example}{Example}

\newcommand{\V}{\operatorname{Var}}
\def\R{\mathbb{R}}
\def\E{\mathbb{E}}
\def\N{\mathbb{N}}
\def\P{\mathbb{P}}
\def\1{\mathbh{1}}
\def\d{{\delta}}
\makeatother

\begin{document}
\begin{frontmatter}

\title{Central limit theorems for $U$-statistics of Poisson point processes}
\runtitle{CLT's for Poisson $U$-statistics}

\begin{aug}
\author[a]{\fnms{Matthias} \snm{Reitzner}\corref{}\ead[label=e1]{matthias.reitzner@uni-osnabrueck.de}}
\and
\author[b]{\fnms{Matthias} \snm{Schulte}\ead[label=e2]{matthias.schulte@uni-osnabrueck.de}}
\runauthor{M. Reitzner and M. Schulte}
\affiliation{University of Osnabrueck and Karlsruher Institut f\"{u}r Technologie}
\address[a]{Institut f\"ur Mathematik\\
Universit\"at Osnabr\"uck\\
49069 Osnabr\"uck\\
Germany\\
\printead{e1}}
\address[b]{Institut f\"ur Stochastik\\
Karlsruher Institut f\"{u}r Technologie\\
76133 Karlsruhe\\
Germany\\
\printead{e2}}
\end{aug}

\received{\smonth{5} \syear{2011}}
\revised{\smonth{7} \syear{2012}}

%
\begin{abstract}
A $U$-statistic of a Poisson point process is defined as the sum $\sum f(x_1, \ldots, x_k)$
over all (possibly infinitely many) $k$-tuples of
distinct points of the point process.
Using the Malliavin calculus, the Wiener--It\^o chaos expansion of such
a functional is computed and used to derive a formula for the variance.
Central limit theorems for $U$-statistics of Poisson point processes are shown,
with explicit bounds for the Wasserstein distance to a Gaussian random variable.
As applications, the intersection process of Poisson hyperplanes and
the length of a random geometric graph are investigated.
\end{abstract}

%
\begin{keyword}[class=AMS]
\kwd[Primary ]{60H07}
\kwd{60F05}
\kwd[; secondary ]{60G55}
\kwd{60D05}
\end{keyword}

\begin{keyword}
\kwd{Central limit theorem}
\kwd{Malliavin calculus}
\kwd{Poisson point process}
\kwd{Stein's method}
\kwd{$U$-statistic}
\kwd{Wiener--It\^o chaos expansion}
\end{keyword}

\end{frontmatter}

\section{Introduction}\label{sec1}
In recent years, Malliavin calculus, Wiener--It\^o chaos expansions and
Fock space representations of functionals of Poisson point processes
have been a rapidly developing topic.
First results already appeared in the classical works of It\^ o \cite
{Ito51,Ito56} and Wiener \cite{Wiener}.
Yet only in the last years prominent contributions produced a deep
theory which most probably will have a strong impact on modern theory
and applications of Poisson point processes; see, for example, Houdre
and Perez-Abreu \cite{HoudrePer}, Last and Penrose \cite
{LastPenrose2011}, Nualart and Vives \cite{NualartVives1990} and Wu
\cite{Wu}.
Here in particular we want to point out the groundbreaking paper by
Peccati et al. \cite{Peccatietal2010} on central limit theorems using
Stein's method and Malliavin calculus. These methods were combined the
first time by Nourdin and Peccati \cite{NourdinPeccati2009} for
functionals depending on Gaussian processes instead of Poisson point
processes. Further developments include the book of Peccati and Taqqu
\cite{PeccatiTaqqu2011} about product formulas for multiple Wiener--It\^
o integrals in the Gaussian and Poisson case and a central limit
theorem due to Peccati and Zheng \cite{PeccatiZheng} generalizing the
main result of \cite{Peccatietal2010} to random vectors.

Poisson point processes occur in many branches of probability theory,
for example, in the theory of Levy processes, and in the theory of
random graphs, in spatial statistics, in communication theory and in
stochastic geometry. Hence there is a wide range of potential
applications of these new results. In this work, we use the Wiener--It\^
o chaos expansion and a related result from \cite{Peccatietal2010} to
prove central limit theorems for a broad class of functionals, namely
for $U$-statistics of Poisson point processes.

Let $\eta$ be a Poisson point process over a state space $X$. We call a
random variable $F $ a $U$-statistic of $\eta$ if
\[
F (\eta) = \sum_{(x_1,\ldots,x_k)\in\eta_{\neq}^k } f(x_1,
\ldots,x_k).
\]
By $\eta_{\neq}^k$ we denote the set of all $k$-tuples of distinct
points of the process.
One should compare this definition to classical $U$-statistics defined on
a set of $n$ random variables $\{ Z_1,\ldots,Z_n\} =\zeta$ where
$ U (\zeta) = \sum_{\zeta_{\neq}^k} f(x_1,\ldots,x_k)$.
For details on classical $U$-statistics we refer to \cite
{Hoeffding1948,KoroljukBorovskich,Lee1990}.
From now on, we mean by $U$-statistic a $U$-statistic of a Poisson point process.

The first step in this paper is the explicit evaluation of expressions
involving Malliavin operators acting on $U$-statistics of Poisson point processes.
The main result of this paper is Theorem \ref{thm:CLT} which gives an
explicit bound on the Wasserstein distance between a normalized
$U$-statistic and a standard Gaussian random variable $N$,
\[
d_W \biggl(\frac{F-\E F}{\sqrt{\V F}},N \biggr)\leq 2k^{7 /2} \sum
_{1\leq i \leq j\leq k}\frac{\sqrt{M_{ij}(f)}}{\V F},
\]
where $M_{ij}(f)$ are sums of certain fourth moment integrals.
If the intensity measure of $\eta$ is of the form $\mu=\lambda\theta$
with an intensity parameter $\lambda\geq1$ and a measure $\theta$, one
is interested in the behavior of $F$ for increasing $\lambda$.
In the particular situation that $f\dvtx X^k\to\R$ is independent of $\lambda
$, we conclude in Theorem \ref{thm:CLTgeometric} that
\[
d_W \biggl(\frac{F-\E F}{\sqrt{\V F}},N \biggr) \leq C_f
\lambda ^{-1/ 2}.
\]
In general this is the optimal rate in $\lambda$ because for a set
$A\subset X$ with $\theta(A)=1$ the $U$-statistic $F=\sum_{x\in\eta} \1(x
\in A)$ is Poisson distributed with parameter $\lambda$, and it is
widely known that a Poisson distributed random variable has this rate
of convergence.

As an application of our result we investigate the intrinsic volumes of
the intersection process of Poisson hyperplanes in a compact convex
window. A~central limit theorem for some of these functionals was
proved in two long and intricate papers by Heinrich \cite{Heinrich2009}
and Heinrich, Schmidt and Schmidt \cite{HeinrichSchmidtSchmidt2006}.
Here we obtain a general result which in addition gives rates of
convergence to Gaussian variables. A second example concerns
functionals of Sylvester type by which we mean the question about the
probability that $k$ points in a convex set are in convex position. Our
last example is about the number of edges of a random geometric graph
in a bounded window. Again we obtain a central limit theorem with a
rate of convergence.
As general references to stochastic geometry and random graphs we refer
to \cite{PenroseRandomGeometricGraphs,SchnWe3} and \cite{SKM}.

To prove our central limit theorems, we first use a result of Last and
Penrose~\cite{LastPenrose2011}, to expand a $U$-statistic in a
Wiener--It\^o chaos expansion as a finite sum of multiple Wiener--It\^o
integrals. This enables us to give a formula for the variance of a
$U$-statistic and to compute two operators from Malliavin calculus that
are defined by their chaos expansions.
Using a theorem for the normal approximation of Poisson functionals due
to Peccati et al. \cite{Peccatietal2010}, we show convergence in the
Wasserstein distance.
In order to apply their result, we need to compute expected values of
products of multiple Wiener--It\^o integrals which is well known to be
a notorious difficult task.
We expect that the same techniques can be used to show central limit
theorems for more general functionals of Poisson point processes.

This paper is organized in the following way. In Section~\ref{sec:Preliminaries}, we introduce Wiener--It\^o chaos expansions for
functionals of a Poisson point process and some operators from
Malliavin calculus. Then we compute the Wiener--It\^o chaos expansion
of a $U$-statistic and its variance in Section~\ref{sec:ChaosExpansionUstatistics}. Using Mallia\-vin calculus we prove
the general version of our central limit theorem for $U$-statistics in
Section~\ref{sec:CLTUstatistics}. Finally, we investigate two special
classes of $U$-statistics and present examples in the Sections~\ref{sec:GeometricUstatistic} and \ref{sec:LocalUstatistic}.


\section{Wiener--It\^o chaos expansions for Poisson point processes}\label{sec:Preliminaries}


\subsection{Poisson point process}
In this paper, we let $\eta$ be a Poisson point process on the measure
space $(X,\mathcal{B}(X),\mu)$ where $X$ is a Borel space and $\mu$ is
a $\sigma$-finite nonatomic Borel measure. A Borel space is a
measurable space which is isomorphic to a Borel subset of $[0,1]$; see
page 7 in \cite{Kallenberg2002}.

More precisely, let $(\Omega,\mathcal{F},P)$ be a probability space.
Denote by $N(X)$ the set of all integer-valued $\sigma$-finite measures
$\nu$ on $X$, equipped with the smallest $\sigma$-algebra $\mathcal
{N}(X)$ such that the mappings $\nu\to\nu(A)$ are measurable for all
sets $A\in\mathcal{B}(X)$.
A random measure $\eta\dvtx  \Omega\to N(X)$ is called a Poisson point
process with intensity measure $\mu$ if for $A\in\mathcal{B}(X)$ the
random variable $\eta(A)$ is Poisson distributed with parameter $\mu
(A)$, and the random variables $\eta(A_1),\ldots,\eta(A_m)$ are
independent for pairwise disjoint sets $A_1,\ldots,A_m\in\mathcal{B}(X)$.
Since the intensity measure $\mu$ is nonatomic, the Poisson point
process is simple, that is, $\eta(\{x\})\leq1$ for all $x\in X$ almost
surely. Thus, we can view $\eta$ as a random set of points in $X$.

As usual, $L^p(X^k)$ denotes the space of all measurable functions $f\dvtx X^k\to\overline{\R}:=\R\cup \{\pm\infty \}$ with
\[
\int_{X^k}\bigl|f(x_1,\ldots,x_k)\bigr|^p
\,d\mu(x_1, \ldots,x_k)<\infty,
\]
where $d\mu(x_1, \ldots,x_k)$ stands for $ d\mu(x_1)\cdots \,d\mu(x_k)$.
Let $L^p_s(X^k)$ be the subset of $\mu^k$-almost everywhere symmetric
functions in $L^p(X^k)$. We call a function symmetric if it is
invariant under all permutations of its arguments. We denote by $\|\cdot
\|$ the norm in $L^2(X^k)$, and by $\langle\cdot,\cdot\rangle$ the
inner product in $L^2(X^k)$. Equipped with this inner product,
$L^2(X^k)$ and $L^2_s(X^k)$ form Hilbert spaces. Instead of the
original probability measure $P$, we always use the image measure $\P=
P \circ\eta$. In the following, $L^p(\P)$ stands for the set of all
measurable functions $F\dvtx N(X)\to\overline{\R}$ with $\E|F|^p<\infty$.

An important property of Poisson point processes is the Slivnyak--Mecke
formula (see Corollary 3.2.3 in \cite{SchnWe3}) which says that
%
%
\begin{equation}
\label{eqn:SM} \E\sum_{(x_1,\ldots,x_k)\in\eta^k_{\neq}}f(x_1,
\ldots,x_k)=\int_{X^k} f(x_1,
\ldots,x_k) \,d\mu(x_1,\ldots,x_k)
\end{equation}
for $f\in L^1(X^k)$. (Recall the definition of $\eta^k_{\neq}$ in the
\hyperref[sec1]{Introduction}.)
The sum on the left-hand side is a priori defined as an $L^1(\P)$ limit
summing only over points in an increasing window. Yet it follows from
the Slivnyak--Mecke formula that $f\in L^1(X^k)$ implies that the sum
on the left-hand side is absolutely convergent almost surely.


\subsection{Multiple Wiener--It\^o integrals}

Now we present the definition of multiple Wiener--It\^o integrals of
order $k\in\N$ following \cite{Surgailis1984}. One starts with simple
functions and extends the definition to arbitrary functions in
$L_s^2(X^k)$. A~function $f\in L^2(X^k)$ is called simple if:
\begin{longlist}[(1)]
\item[(1)] $f$ is symmetric;
\item[(2)] $f$ is constant on a finite number of Cartesian products
$B_1\times\cdots\times B_k\in\mathcal{B}(X)^k$ and vanishes elsewhere;
\item[(3)] $f$ vanishes on diagonals, which means $f(x_1,\ldots,x_k)=0$
if $x_i= x_j$ for some $i\neq j$.
\end{longlist}
Let $\mathcal{S}(X^k)$ be the space of all simple functions. For
$f_0\in\mathcal{S}(X^k)$ and $k \in\N$, the multiple Wiener--It\^o
integral $I_k(f_0)$ of $f_0$ with respect to the compensated Poisson
point process $\eta-\mu$ is defined by
\[
I_k(f_0)=\int_{X^k} f_0
\,d(\eta-\mu)^k =\sum f_0^{B_1\times
\cdots\times B_k}(\eta-\mu)
(B_1)\cdots(\eta-\mu) (B_k),
\]
where we sum over all Cartesian products and $f_0^{B_1\times\cdots
\times B_k}$ is the value of $f_0$ on such a set. For $k=0$ we put
$I_0(f)=f$. By a straightforward computation, one shows
%
%
\begin{equation}
\E I_k(f_0)^2=k!\|f_0
\|^2.\label{eqn:ISO}
\end{equation}
Thus there is an isometry between $\mathcal{S}(X^k)$ and a subset of
$L^2(\P)$. Furthermore, $\mathcal{S}(X^k)$ is dense in $L^2_s(X^k)$,
whence for every $f\in L^2_s(X^k)$ there is a sequence\vadjust{\goodbreak} $(f_n)_{n\in
\mathbb{N}}$ of simple functions with $f_n\rightarrow f$ in
$L^2_s(X^k)$. Because of the isometry (\ref{eqn:ISO}), it is possible
to define $I_k(f)$ as the limit of $(I_k(f_n))_{n\in\mathbb{N}}$ in
$L^2(\P)$. Hence for an arbitrary symmetric function $f \in L^2_s
(X^k)$ we put $f^0(x_1,\ldots,x_k)=f(x_1,\ldots,x_k)$ if $x_i\neq x_j$
for all $i\neq j$ and $f^0(x_1,\ldots,x_k)=0$ otherwise and obtain
\[
I_k(f)=\int_{X^k} f^{0} \,d(\eta-
\mu)^k.
\]

We remark that the denseness of $\mathcal{S}(X^k)$ in $L_s^2(X^k)$
depends on the topological structure of $X$ and the fact that $\mu$ is
nonatomic. For a definition without these requirements we refer to
\cite{LastPenrose2011}.

It follows directly from the definition that multiple Wiener--It\^o
integrals have the properties summarized in the following:
%
\begin{lemma}\label{lem:ItoIntegrals}
Let $f\in L^2_s(X^n)$ and $g\in L^2_s(X^m)$ with $n,m\geq1$. Then:
\begin{longlist}[(a)]
\item[(a)] $ \E I_n(f)=0$;
\item[(b)] $ \E I_n(f)I_m(g)=\1(n = m) n ! \langle f, g \rangle$.
\end{longlist}
\end{lemma}


\subsection{Wiener--It\^o chaos expansions}\label{sec:WIchaos}

For a measurable function $F\dvtx N(X)\to\overline{\R}$ and $y\in X$ we
define the difference operator as
\[
D_yF(\eta)=F(\eta+\delta_y)-F(\eta),
\]
where $\delta_y$ is the Dirac measure at the point $y$. The difference
operator $D_yF$ measures the effect of adding the point $y\in X$ to the
Poisson point process, whence it is also denoted as add one cost
operator in \cite{LastPenrose2011}. The iterated difference operator is
defined by
\[
D_{y_1,\ldots,y_i}F=D_{y_1}D_{y_2,\ldots,y_i}F.
\]
Let the functions $f_i\dvtx X^i\to\overline{\R}$ be given by $f_0=\E F$ and
\[
f_i(y_1,\ldots,y_i)=\frac{1}{i!}\E
D_{y_1,\ldots,y_i}F,\qquad   i \geq1,
\]
if these expectations exist.
Because of the symmetry of the iterated difference operator, $f_i$ is
symmetric if defined. The following relationships between $F$, the
functions $f_i$, $i\in\N$, and the variance of $F$ have been shown by
Last and Penrose \cite{LastPenrose2011}.

\begin{theorem}[(Last and Penrose \cite{LastPenrose2011})]\label{prop:LastPenrose}
Let $F\in L^2(\P)$. Then $f_i\in L_s^2(X^i)$, $i\in\mathbb{N}$ and
\[
F=\sum_{i=0}^{\infty}I_i(f_i),
\]
where the sum converges in $L^2(\P)$. The $f_i\in L_s^2(X^i),i\in\N$
are the $\mu^i$-almost everywhere unique $g_i\in L_s^2(X^i)$, $i\in\N$,
satisfying\vadjust{\goodbreak} $F=\sum_{i=0}^\infty I_i(g_i)$ in $L^2(\P)$. Furthermore,
\[
\V F= \sum_{i=1}^\infty i!\|f_i
\|^2.
\]
\end{theorem}
In the following, we call the functions $f_i$, $i\in\N$, kernels of the
Wiener--It\^o chaos expansion of $F$. The class of sequences $
(g_i )_{i\in\N}$ with $g_i\in L_s^2(X^i)$ and
\[
\sum_{i=0}^\infty i! \|g_i
\|^2<\infty
\]
composes a Hilbert space isomorphic to the symmetric Fock space
associated with $L^2(X)$. In this context, Theorem \ref
{prop:LastPenrose} states that there exists an isometry between $L^2(\P
)$ and a symmetric Fock space.


\subsection{Malliavin calculus}
Our proofs for central limit theorems are based on a result for the
normal approximation of Poisson functionals from \cite
{Peccatietal2010}, which uses operators from Malliavin calculus. In the
following, we give a short introduction to these operators. For more
details we refer to \cite{LastPenrose2011,NualartVives1990,Peccatietal2010}.

Let $F\in L^2(\P)$ and $f_i$, $i\in\N$, be the kernels of the
Wiener--It\^o chaos expansion of $F$. First of all, we give an
alternative definition of the difference operator $D_y$ using the
Wiener--It\^o chaos expansion of $F$.
%
\begin{definition}
Let
%
%
\begin{equation}
\label{eqn:ConditionDifference} \sum_{i=1}^\infty
i i!\|f_i\|^2<\infty.
\end{equation}
Then the random function $y\mapsto D_yF, y\in X$, is given by
\[
D_yF=\sum_{i=1}^\infty i
I_{i-1}\bigl(f_i(y,\cdot)\bigr).
\]
\end{definition}
It can be proved (see \cite{NualartVives1990}, Theorem 6.2 or
\cite{LastPenrose2011}, Theorem 3.3) that for $F\in L^2(\P)$ satisfying (\ref
{eqn:ConditionDifference}) this definition coincides with the one
introduced in Section~\ref{sec:WIchaos}.

\begin{definition}
If
\[
\sum_{i=1}^\infty i^2i!
\|f_i\|^2<\infty,
\]
then the Ornstein--Uhlenbeck generator $LF$ is the random variable
given by
\[
LF=-\sum_{i=1}^{\infty}i I_i(f_i).
\]
\end{definition}
The Ornstein--Uhlenbeck generator has an inverse operator. Its domain
is the space of all centred $F\in L^2(\P)$, that is, $F\in
L^2(\P)$ with $\E F=0$, and
\[
L^{-1}F=-\sum_{i=1}^\infty
\frac{1}{i}I_i(f_i).
\]

If $F$ is in the domain of $L$, then the Ornstein--Uhlenbeck generator
can be written as
%
%
\begin{equation}
\label{eqn:defL} LF = \int_X F( \eta-
\delta_x) - F(\eta) \,d\eta(x) - \int_X \bigl(F(
\eta) - F(\eta+ \delta_z)\bigr) \,d \mu(z).
\end{equation}
This follows from the representation of the difference operator and the
Skorohod-integral (see \cite{LastPenrose2011}, formula (3.19)), which
is not used in this work.


\section{Malliavin calculus and Wiener--It\^o chaos expansions for
$U$-statistics}\label{sec:ChaosExpansionUstatistics}

In this section, we define $U$-statistics of Poisson point processes and
investigate their Wiener--It\^o chaos expansions.
In particular, we apply the Malliavin operators to $U$-statistics and
present explicit formulae for the kernels of the Wiener--It\^o chaos
expansion and the variance.

\subsection{$U$-statistics of Poisson point processes}

Recall the definition $\eta_{\neq}^k =  \{(x_1,\ldots,x_k)\in\eta
^k, x_i\neq x_j \mbox{ for }i\neq j \}$ from the \hyperref[sec1]{Introduction}.

\begin{definition}\label{def:DefinitionUStatistic}
A random variable
%
%
\begin{equation}
\label{eqn:DefinitionUStatistic} F=\sum_{(x_1,\ldots,x_k)\in\eta^k_{\neq}}f(x_1,
\ldots,x_k)
\end{equation}
with $f\in L_s^1(X^k)$ is called $U$-statistic of order $k$.
\end{definition}

By the Slivnyak--Mecke formula (\ref{eqn:SM}), it holds that
\[
\E\sum_{(x_1,\ldots,x_k)\in\eta^k_{\neq}}f(x_1,\ldots,x_k)=
\int_X\cdots\int_X
f(x_1,\ldots,x_k) \,d\mu(x_1, \ldots,
x_k)
\]
so that $f\in L_s^1(X^k)$ guarantees $F\in L^1(\P)$.
Due to the fact that we sum over all permutations of $k$ points in (\ref
{eqn:DefinitionUStatistic}), we can assume without loss of generality
in Definition~\ref{def:DefinitionUStatistic} that $f$ is symmetric.

Since we want to use Wiener--It\^o chaos expansions, we always require
that $F$ is in $L^2(\P)$. For the central limit theorems we
additionally assume that $F$ is absolutely convergent.
%
\begin{definition}\label{def:absolutelyconvergent}
A $U$-statistic $F $ is absolutely convergent if
\[
\overline{F}=\sum_{(x_1,\ldots,x_k)\in\eta^{k}_{\neq}}\bigl|f(x_1,
\ldots,x_k)\bigr|
\]
is in $L^2(\P)$.\vadjust{\goodbreak}
\end{definition}
Note that $F$ absolutely convergent implies that $F \in L^2(\P)$.
Obviously every $F\in L^2(\P)$ with $f\geq0$ is absolutely convergent.

\subsection{Malliavin calculus}
We start by calculating the difference operator of a $U$-statistic $F$.
%
\begin{lemma}\label{Lemma:Dy}
Let $F\in L^2(\P)$ be a $U$-statistic of order $k$. Then the difference
operator applied to $F$ gives
\[
D_{y_1}F = k\sum_{(x_1,\ldots,x_{k-1})\in\eta_{\neq
}^{k-1}}f(y_1,x_1,
\ldots,x_{k-1}).
\]
\end{lemma}
\begin{pf}
By the definition of the difference operator $D_y$ and the symmetry of~$f$, we obtain for a $U$-statistic
\begin{eqnarray*}
D_{y_1}F &=& \sum_{(x_1,\ldots,x_k)\in(\eta\cup \{y_1 \})_{\neq
}^k}f(x_1,
\ldots,x_k)-\sum_{(x_1,\ldots,x_k)\in\eta_{\neq
}^k}f(x_1,
\ldots,x_k)
\\
&=& \sum_{(x_1,\ldots,x_{k-1})\in\eta_{\neq}^{k-1}} \bigl(f(y_1,x_1,
\ldots,x_{k-1})+\cdots+f(x_1,\ldots,x_{k-1},y_1)
\bigr)
\\
& =& k\sum_{(x_1,\ldots,x_{k-1})\in\eta_{\neq}^{k-1}}f(y_1,x_1,
\ldots,x_{k-1}).
\end{eqnarray*}
\upqed\end{pf}

An analogous straightforward computation using (\ref{eqn:defL})
verifies the following lemma.
%
\begin{lemma}\label{Lemma:LF}
Let $F\in L^2(\P)$ be a $U$-statistic of order $k$. Then the
Ornstein--Uhlenbeck operator applied to $F$ gives
\[
LF = -kF + k \int_X \sum_{(x_1,\ldots,x_{k-1})\in\eta
_{\neq}^{k-1}}
f(x_1, \ldots, x_{k-1}, z) \,d \mu(z).
\]
\end{lemma}

Without proof we also state the inverse Ornstein--Uhlenbeck operator of
a $U$-statistic.
\begin{eqnarray*}
&&L^{-1}  (F- \E F)\\
&&\qquad= \Biggl(\sum_{m=1}^k
\frac1m \Biggr) \int_{X^k} f(y_1,
\ldots,y_{k}) \,d\mu(y_1, \ldots, y_k)
\\
& &\qquad\quad{}  - \sum_{m=1}^{k} \frac1m
\sum_{(x_1,\ldots,x_{m})\in\eta_{\neq}^{m}} \int_{X^{k-m}}
f(x_{1},\ldots,x_m,y_1,\ldots,y_{k-m})
\\
&&\hspace*{160pt}d\mu(y_1, \ldots, y_{k-m }).
\end{eqnarray*}

\subsection{Wiener--It\^o chaos expansions}
Let us now compute the kernels and the Wiener--It\^o chaos expansion of
a $U$-statistic
$F=\sum_{\eta_{\neq}^k} f$ with $F \in L^2(\P)$.

\begin{lemma}\label{Lemma:Kernels}
Let $F\in L^2(\P)$ be a $U$-statistic of order $k$. Then the kernels of
the Wiener--It\^o chaos expansion of $F$ have the form\vspace*{1pt}
\[
f_i(y_1,\ldots,y_i)= %
\cases{
\displaystyle\pmatrix{k
\cr
i}\int_{X^{k-i}}f(y_1,
\ldots,y_i,x_1,\ldots,x_{k-i}) \,d
\mu(x_1,\ldots,x_{k-i}), \vspace*{2pt}\cr
 \quad\hspace*{22pt} i\leq k,\vspace*{2pt}
\cr
0,\qquad
i>k,} \vspace*{1pt}%
\]
and $F$ has the variance\vspace*{1pt}
%
%
\begin{eqnarray}
\label{eq:GeneralVariance} \V F  &= & \sum_{i=1}^{k}i!
\pmatrix{k
\cr
i}^2
\nonumber\\[1pt]
&&\hspace*{15pt}\times{} \int_{X^i} \biggl( \int
_{X^{k-i}} f(y_1,\ldots,y_i,
x_1,
\ldots,x_{k-i}) \,d\mu(x_1, \ldots, x_{k-i})
\biggr)^2 \\[1pt]
&&\hspace*{43pt}d\mu(y_1, \ldots, y_i).\nonumber\vspace*{1pt}
\end{eqnarray}
\end{lemma}

For the special case $k=2$ the formulas for the kernels are already
implicit in the paper by Molchanov and Zuyev \cite{MolZuy} where ideas
closely related to Malliavin calculus have been used.\vspace*{1pt}

\begin{pf*}{Proof of Lemma \ref{Lemma:Kernels}}
In Lemma \ref{Lemma:Dy}, the difference operator of a $U$-statistic was computed.
Proceeding by induction, we get
\[
D_{y_1,\ldots,y_{i}}F  =  \frac{k!}{(k-i)!}\sum_{(x_1,\ldots,x_{k-i})\in\eta_{\neq
}^{k-i}}f(y_{1},
\ldots,y_{i},x_1,\ldots,x_{k-i})\vspace*{1pt}
\]
for $i\leq k$. Hence $D_{y_1,\ldots,y_k}F$ only depends on $y_1,\ldots,y_k$
and is independent of the Poisson point process. This yields\vspace*{1pt}
\[
D_{y_1,\ldots,y_{k+1}}F =0\quad \mbox{and}\quad  D_{y_1,\ldots,y_i}F=0\vspace*{1pt}
\]
for all $i>k$.
We just proved\vspace*{1pt}
\begin{eqnarray*}
&&D_{y_1,\ldots,y_i}F\\
&&\qquad= %
\cases{\displaystyle\frac{k!}{(k-i)!}\sum
_{(x_1,\ldots,x_{k-i})\in\eta_{\neq
}^{k-i}}f(y_1,\ldots,y_i,x_1,
\ldots,x_{k-i}), &\quad $i\leq k,$\vspace*{2pt}
\cr
0, &\quad $\mbox{otherwise.}$}
\end{eqnarray*}
By the Slivnyak--Mecke formula (\ref{eqn:SM}), we obtain
\begin{eqnarray*}
f_i(y_1,\ldots,y_i) &=&\frac{1}{i!}\E
D_{y_1,\ldots,y_i}F
\\
&=&\frac
{1}{i!}\E\frac{k!}{(k-i)!}\sum_{(x_1,\ldots,x_{k-i})\in\eta_{\neq
}^{k-i}}f(y_1,
\ldots,y_i,x_1,\ldots,x_{k-i})
\\
&=&\frac{k!}{i!(k-i)!}\int_{X^{k-i}}f(y_1,
\ldots,y_i,x_1,\ldots,x_{k-i}) \,d
\mu(x_1,\ldots,x_{k-i})
\end{eqnarray*}
for $i\leq k$. The formula for the variance follows from Proposition
\ref{prop:LastPenrose}.
\end{pf*}

Note that $F \in L^ 2(\P)$ implies $f_i \in L^2_s (X^ i)$, and thus
that for all $1\leq i \leq k$
\[
\int_{X^{i}} \biggl(\int_{X^{k-i}}f(y_1,
\ldots,y_i,x_1,\ldots,x_{k-i}) \,d
\mu(x_1,\ldots,x_{k-i}) \biggr) ^2 \,d
\mu(y_1,\ldots,y_{i}) < \infty.
\]
In particular, it holds $f\in L_s^2(X^k)$.

By Lemma~\ref{Lemma:Kernels}, $U$-statistics only have a finite number of
nonvanishing kernels.
The following theorem characterizes a $U$-statistic by this property. We
call a Wiener--It\^o chaos expansion finite if only a finite number of
kernels do not vanish.

\begin{theorem}\label{thm:UStatisticIto}
Assume $F\in L^2(\P)$.
\begin{longlist}[(1)]
\item[(1)] If $F$ is a $U$-statistic, then $F$ has a finite Wiener--It\^o
chaos expansion with kernels $f_i\in L_s^1(X^i)\cap L_s^2(X^i) $,
$i=1,\ldots,k$.
\item[(2)] If $F$ has a finite Wiener--It\^o chaos expansion with
kernels $f_i\in L_s^1(X^i) \cap L_s^2 (X^ i)$, $i=1,\ldots,k$, then $F$
is a (finite) sum of $U$-statistics and a constant.
\end{longlist}
\end{theorem}
\begin{pf}
The fact that a $U$-statistic $F\in L^2(\P)$ has a finite Wiener--It\^o
chaos expansion with $f_i\in L_s^1(X^i)$ follows from Lemma \ref
{Lemma:Kernels} and from $f \in L_s^1(X^k)$.

For the second part of the proof, let $F\in L^2(\P)$ have a finite
Wiener--It\^o chaos expansion, that is,
\[
F=\sum_{i=0}^m I_i(f_i)
\]
with kernels $f_i\in L^1_s(X^i)\cap L^2_s(X^i)$ and $m\in\N$. Now
Proposition 4.1 in \cite{Surgailis1984} implies that
\[
I_i(f_i)=\sum_{j=0}^i
(-1)^{i-j} \pmatrix {i
\cr
j} \sum_{(x_1,\ldots,x_j)\in
\eta^j_{\neq}}
f_i^{(j)}(x_1,\ldots,x_j),
\]
where the inner sum is a constant for $j=0$ and $f_i^{(j)}$ is given by
\[
f_i^{(j)}(x_1,\ldots,x_j)=\int
_{X^{i-j}} f_i(x_1,\ldots,x_j,y_1,\ldots,y_{i-j}) \,d\mu(y_1,
\ldots,y_{i-j}).
\]
The assumption $f_i\in L^1_s(X^i)$ guarantees $f_i^{(j)}\in L^1_s(X^j)$
for $j=1,\ldots,i$ and \mbox{$f_i^{(0)}\in\R$}. Hence, every Wiener--It\^o
integral is a (finite) sum of $U$-statistics and a constant, and the same
holds for $F$.
\end{pf}

\subsection{Examples} The following examples show that the assumptions
on $F$ and~$f_i$ in Theorem \ref{thm:UStatisticIto} are necessary. In
all examples, we consider a Poisson point process in $\R$ with the
Lebesgue measure as intensity measure.

\begin{example*}
There exist random variables in $L^ 2(\P)$ with finite Wiener--It\^o
chaos expansions which are not sums of $U$-statistics. This is possible
if the kernels $f_i$ are in $L_s^ 2(X^ i) \setminus L_s^ 1 (X^ i)$.
Define $g\dvtx \R\to\R$ as
\[
g(x)= \frac{1}{x} \1\bigl( |x|>1\bigr),
\]
which is in $ L^2(\R) \setminus L^1(\R)$. Now we define the random
variable $G=I_1(g)$. $G$ is in $L^2(\P)$ and has a finite Wiener--It\^o
chaos expansion. But the formal representation
\[
I_1(g)=\sum_{x\in\eta} g(x) -\int
_{\R}g(x) \,dx
\]
we used in the proof of Theorem \ref{thm:UStatisticIto} fails because
the integral does not exist.
\end{example*}

\begin{example*}
There also exist $U$-statistics $F\in L^1(\P)$ with $f\in L_s^1(X^k)\cap
L_s^2(X^k)$ which are not in $L^2(\P)$.
We construct $f\in L^1_s(\R^2)\cap L^2_s(\R^2)$ with $\| f_1\|= \infty$
by putting
\[
f(x_1,x_2)= \1 (0 \leq x_1 \sqrt
x_2 \leq1) \1 (0 \leq x_2 \sqrt x_1 \leq1)
\]
and
\[
F=\sum_{(x_1,x_2)\in\eta_{\neq}^2}f(x_1,x_2).
\]
In this case the first kernel,
\[
f_1(y)=\E\biggl[2\sum_{x\in\eta} f(y,x)
\biggr]=2\int_{\R} f(y,x) \,dx=2\1 (y\geq0)\min \biggl\{
\frac{1}{y^2},\frac{1}{\sqrt{y}} \biggr\}
\]
is not in $L_s^2(\R)$ so that $F$ has no Wiener--It\^o chaos expansion
and cannot be in $L^2(\P)$.
\end{example*}

\begin{example*}
By Theorem \ref{thm:UStatisticIto}(2), a functional $F \in L^ 2(\P)$
with a finite Wiener--It\^o chaos expansion and kernels $f_i\in
L_s^1(X^i) \cap L_s^ 2 (X^ i)$, $i=1,\ldots,k$, is a (finite) sum of
$U$-statistics. Our next example shows that neither the single
$U$-statistics are in $L^2(\P)$ nor are the summands necessarily in $L^2_s(X^i)$.
Set $F=I_2(f)$ with $f$ as above. Then
\[
I_2(f)=\int_{\R^2 } f(x,y) \,dx \,dy- 2\sum
_{x\in\eta}\int_{\R} f(x,y) \,dy+\sum
_{(x_1,x_2)\in\eta
^2_{\neq}}f(x_1,x_2),
\]
and $F$ is a sum of $U$-statistics. Since
$\E[(\sum_{x\in\eta}\int_{\R} f(x,y)  \,dy)^2]=\infty$, we know
that the $U$-statistic
\[
\sum_{x\in\eta}\int_{\R} f(x,y) \,dy
\]
is not in $L^2(\P)$, nor are the summands $\int_{\R} f(x,y)  \,dy$ in
$L^2(\R)$. This is in contrast to the remark after the proof of Lemma
\ref{Lemma:Kernels} that for a $U$-statistic $F \in L^2(\P)$, we always
have $f \in L^2(X^k)$.
\end{example*}

\begin{example*}
To motivate the definition of an absolutely convergent $U$-statistic, we
give an example of a $U$-statistic that is in $L^2(\P)$ but not
absolutely convergent.
Similarly to the previous examples, we set
\[
f(x_1,x_2)=\1\bigl(0\leq|x_1|
\sqrt{|x_2|}\leq1\bigr)\1\bigl(0\leq|x_2|\sqrt {|x_1|}
\leq1\bigr) \bigl(2\1(x_1x_2\geq0)-1 \bigr)
\]
and
\[
F=\sum_{(x_1,x_2)\in\eta^2_{\neq}}f(x_1,x_2)
\quad\mbox{and}\quad \overline{F}=\sum_{(x_1,x_2)\in\eta^2_{\neq}}\bigl|f(x_1,x_2)\bigr|.
\]
Now it is easy to verify that $f_1(x)=0$ and $f_2(x_1,x_2)=f(x_1,x_2)$
so that $F\in L^2(\P)$. But the first kernel of $\overline{F}$ is not
in $L^2(\R)$ so that $\overline{F}\notin L^2(\P)$.
\end{example*}


\section{Central limit theorems for $U$-statistics}\label{sec:CLTUstatistics}

In this section, we derive a central limit theorem for $U$-statistics of
Poisson point processes. In particular, we are interested in the
Wasserstein distance of a normalized $U$-statistic and a standard
Gaussian random variable. Recall that the Wasserstein distance
$d_W(Y,Z)$ of two random variables $Y$ and $Z$ is given by
\[
d_W(Y,Z)=\sup_{h\in \operatorname{Lip}(1)}\bigl|\E h(Y)-\E h(Z)\bigr|,
\]
where $\operatorname{Lip}(1)$ is the set of all functions $h\dvtx \R\to\R$ with a
Lipschitz-constant less than or equal to one.
It is important to note that convergence in the Wasserstein distance
implies convergence in distribution. In particular, it is known (see
\cite{ChenShao2005}, e.g.) that for a Gaussian random variable $N$ we have
\[
\bigl|\P(Y\leq t)-\P(N\leq t)\bigr| \leq2 \sqrt{d_{W}(Y,N)}
\]
for all $t\in\R$. Hence, we can prove central limit theorems by showing
convergence to a Gaussian random variable in the Wasserstein distance.\vadjust{\goodbreak}

Our main estimate for the distance between $F= \sum_{\eta^k_{\neq}} f$
and a standard Gaussian random variable $N$ is Theorem~\ref{thm:CLT}
which states that
\[
d_W \biggl(\frac{F-\E F}{\sqrt{\V F}},N \biggr)\leq 2k^{7 /2} \sum
_{1\leq i \leq j\leq k}\frac{\sqrt{M_{ij}(f)}}{\V
F},
\]
where the $M_{ij} (f) $ are sums of certain fourth moment integrals.
The precise definition is given in formula (\ref{def:Mij}). In most
applications, it is elementary to bound these fourth moments of $f$.
This is carried out in Sections~\ref{sec:GeometricUstatistic} and~\ref{sec:LocalUstatistic}.


\subsection{An abstract CLT}

Our most general result is the following upper bound for the
Wasserstein distance of a Poisson functional with a finite Wiener--It\^
o chaos expansion and a standard Gaussian random variable.
To neatly formulate our results and proofs, we use the abbreviations
%
%
\begin{eqnarray}\label{eqn:R1}
R_{ij} & = & \E \biggl(\int_{X}
I_{i-1}\bigl(f_i(z,\cdot)\bigr) I_{j-1}
\bigl(f_j(z,\cdot)\bigr) \,d\mu(z) \biggr)^2
\nonumber
\\[-8pt]
\\[-8pt]
\nonumber
& &{} - \biggl[ \E\int_X I_{i-1}
\bigl(f_i(z,\cdot)\bigr) I_{j-1}\bigl(f_j(z,
\cdot)\bigr) \,d\mu(z) \biggr]^2,
\\
\tilde{R}_{i} & = & \E\int_X I_{i-1}
\bigl(f_i(z,\cdot)\bigr)^4 \,d\mu(z) \label{eqn:R2}
\end{eqnarray}
for $i,j=1,\ldots,k$. Note that $R_{11}=0$ and that for $i \neq j$ the
second expectation in $R_{ij}$ vanishes.

\begin{theorem}\label{thm:GeneralCLT}
Suppose $F \in L^ 2(\P)$ has a finite Wiener--It\^o chaos expansion of
order $k$, and $N$ is a standard Gaussian random variable. Then
%
%
\begin{equation}
\label{eqn:thgenCLT} \,d_W \biggl(\frac{F-\E F}{\sqrt{\V F}},N \biggr) \leq k
\sum_{1\leq i, j\leq k}\frac{\sqrt{R_{ij}}}{\V F}+k^{7 /2}\sum
_{i=1}^k\frac{\sqrt{\tilde{R}_i}}{\V F}
\end{equation}
with $R_{ij}$ and $\tilde{R}_i$ defined in (\ref{eqn:R1}) and (\ref{eqn:R2}).
\end{theorem}

\begin{pf}
Our proof is based on the following result of Peccati et al.
(Theorem~3.1 in \cite{Peccatietal2010}), which is derived by a
combination of Malliavin calculus and Stein's method.
%
\begin{theorem}[(Peccati et al. \cite{Peccatietal2010})]\label{thm:CLTPeccati}
Let $G\in L^2(\P)$ with $\E G=0$ be in the domain of $D$ and let $N$ be
a standard Gaussian random variable. Then
\begin{eqnarray*}
d_W(G,N) & \leq& \E\bigl|1-\bigl\langle DG,-DL^{-1}G\bigr
\rangle\bigr|+\int_X \E \bigl[|D_zG|^2\bigl|D_zL^{-1}G\bigr|
\bigr] \,d\mu(z)
\\
&\leq& \sqrt{\E \bigl(1-\bigl\langle DG,-DL^{-1}G\bigr\rangle
\bigr)^2}+\int_X \E \bigl[|D_zG|^2\bigl|D_zL^{-1}G\bigr|
\bigr] \,d\mu(z).
\end{eqnarray*}
\end{theorem}
From now on, we denote by
\[
G=\frac{F-\E F}{\sqrt{\V F}}
\]
the normalization of $F$ and by $g_i\in L^2_s(X^i),i=1,\ldots,k$, the
kernels of $G$. Thus, it follows
\[
g_i(x_1,\ldots,x_i)=\frac{1}{\sqrt{\V F}}f_i(x_1,
\ldots,x_i)
\]
for $i=1,\ldots,k$ and
$\V G=\sum_{i=1}^ki!\|g_i\|^2=1$.

Since $F$ has a finite Wiener--It\^o chaos expansion, $F$ is in the
domain of $D$, and we can apply the above theorem. By the definitions
of the Malliavin operators and the triangle inequality, we obtain
\begin{eqnarray*}
&&\E\bigl\llvert 1-\bigl\langle DG, -DL^{-1}G\bigr\rangle\bigr
\rrvert
\\
&&\qquad= \E \Biggl|\sum_{i=1}^k i!
\|g_i\|^2-\int_X \sum
_{i=1}^kiI_{i-1}\bigl(g_i(z,
\cdot)\bigr) \sum_{i=1}^k
I_{i-1}\bigl(g_i(z,\cdot)\bigr) \,d\mu(z) \Biggr|
\\
&&\qquad \leq \sum_{i=2}^k \E \biggl| i!
\|g_i\|^2-i\int_{X}
I_{i-1}\bigl(g_i(z,\cdot)\bigr) I_{i-1}
\bigl(g_i(z,\cdot )\bigr) \,d\mu(z) \biggr|
\\
&&\qquad\quad{} +\sum_{i,j=1,  i \neq j}^k i\E \biggl|\int
_X I_{i-1}\bigl(g_i(z,\cdot )\bigr)
I_{j-1}\bigl(g_j(z,\cdot)\bigr) \,d\mu(z)\biggr |.
\end{eqnarray*}
The first sum on the right-hand side of the inequality starts with
$i=2$ since the summand for $i=1$ vanishes. As a consequence of
Fubini's theorem and Lemma~\ref{lem:ItoIntegrals}, it holds that
\[
\E i\int_{X} I_{i-1}\bigl(g_i(z,
\cdot)\bigr) I_{i-1}\bigl(g_i(z,\cdot)\bigr) \,d\mu(z)=i!
\|g_i\|^2.
\]
Combining this with the Cauchy--Schwarz inequality leads to
\begin{eqnarray*}
&&\E\biggl |i!\|g_i\|^2-i\int_{X}
I_{i-1}\bigl(g_i(z,\cdot )\bigr)I_{i-1}
\bigl(g_i(z,\cdot)\bigr) \,d\mu(z) \biggr|
\\
&&\qquad\leq \sqrt{\E \biggl(i!\|g_i\|^2-i\int
_{X} I_{i-1}\bigl(g_i(z,\cdot )
\bigr)I_{i-1}\bigl(g_i(z,\cdot)\bigr) \,d\mu(z)
\biggr)^2}
\\
&&\qquad= \sqrt{i^2 \E \biggl(\int_{X}
I_{i-1}\bigl(g_i(z,\cdot)\bigr) I_{i-1}
\bigl(g_i(z,\cdot)\bigr) \,d\mu(z) \biggr)^2
-(i!)^2\|g_i\|^4}\\
&&\qquad= i \frac{\sqrt{R_{ii}}}{\V F}
\end{eqnarray*}
and
\begin{eqnarray*}
&&\E\biggl |\int_{X} I_{i-1}\bigl(g_i(z,
\cdot)\bigr)I_{j-1}\bigl(g_j(z,\cdot)\bigr) \,d\mu(z) \biggr|
\\
&&\qquad \leq \sqrt{\E \biggl(\int_X I_{i-1}
\bigl(g_i(z,\cdot)\bigr) I_{j-1}\bigl(g_j(z,
\cdot)\bigr) \,d\mu(z) \biggr)^2} = \frac{\sqrt{R_{ij}}}{\V F}
\end{eqnarray*}
for $i\neq j$. Now it holds that
%
%
\begin{eqnarray}
\label{eqn:term1} \E\bigl\llvert 1-\bigl\langle DG, -DL^{-1}G\bigr
\rangle\bigr\rrvert & \leq& \sum_{i=2}^k i
\frac{\sqrt{R_{ii}}}{\V F} + \sum_{i,j=1, i \neq j}^k i
\frac{\sqrt {R_{ij}}}{\V F}
\nonumber
\\[-8pt]
\\[-8pt]
\nonumber
& \leq& k\sum_{1\leq i,j\leq k}\frac{\sqrt {R_{ij}}}{\V F}.
\nonumber
\end{eqnarray}
Furthermore, again by the Cauchy--Schwarz inequality we have
\begin{eqnarray*}
&&\int_X \E\bigl[(D_zG)^2\bigl|D_zL^{-1}G\bigr|
\bigr] \,d\mu(z)
\\
&&\qquad\leq \biggl( \int_X \E\bigl[(D_zG)^4
\bigr] \,d\mu(z) \biggr)^{{1}/{2}} \biggl( \int_X \E
\bigl[\bigl(D_zL^{-1}G\bigr)^2\bigr] \,d\mu(z)
\biggr)^{{1}/{2}}.
\end{eqnarray*}
By the definitions of the Malliavin operators and H\"older's
inequality, we can rewrite the expressions on the right-hand side as
\begin{eqnarray*}
\int_X \E\bigl[\bigl(D_zL^{-1}G
\bigr)^2\bigr] \,d\mu(z)&=& \int_X\sum
_{i=1}^k \E\bigl[I_{i-1}
\bigl(g_i(z,\cdot)\bigr)^2\bigr] \,d\mu(z)\\
& = &\sum
_{i=1}^{k} (i-1)!\|g_i\|^2
\leq1
\end{eqnarray*}
and
\[
\int_X \E\bigl[(D_zG)^4\bigr] \,d
\mu(z) \leq \int_X k^3 \sum
_{i=1}^k i^4 \E\bigl[I_{i-1}
\bigl(g_i(z,\cdot)\bigr)^4\bigr] \,d\mu(z)=k^3
\sum_{i=1}^k i^4
\frac{\tilde{R}_i}{(\V F)^2}.
\]
Hence
%
%
\begin{eqnarray}
\label{eqn:term2} \int_X\E\bigl[(D_zG)^2\bigl|D_zL^{-1}G\bigr|
\bigr] \,d\mu(z) & \leq& \sqrt{k^3 \sum
_{i=1}^k i^4 \frac{\tilde{R}_i}{(\V F)^2}}
\nonumber
\\[-8pt]
\\[-8pt]
\nonumber
& \leq& k\sqrt{k}\sum_{i=1}^k
i^2\frac{\sqrt{\tilde{R}_i}}{\V F}\leq k^{7 /2}\sum
_{i=1}^k\frac{\sqrt{\tilde{R}_i}}{\V F}.
\nonumber
\end{eqnarray}
Combining Theorem \ref{thm:CLTPeccati} with formulas (\ref{eqn:term1})
and (\ref{eqn:term2}) gives the right-hand side of~\eqref{eqn:thgenCLT}
in Theorem \ref{thm:GeneralCLT}.
\end{pf}


\subsection{Estimates for the error terms}
To estimate the right-hand side of (\ref{eqn:thgenCLT}) for a
$U$-statistic\vspace*{-1pt} $F=\sum_{\eta^k_{\neq}}f$ in terms of the function $f$, we
are interested in the behavior of $R_{ij}$ and $\tilde{R}_i$ for
$i,j=1,\ldots,k$. Thus we need to compute expected values of the type
\[
\E\prod_{l=1}^m I_{n_l}(f_l),\qquad
 m\in\N,  n_1,\ldots,n_m\in\N
\]
with $f_l\in L_s^1(X^{n_l})\cap L_s^2(X^{n_l})$ for $l=1,\ldots,m$.
Such products of multiple Wiener--It\^o integrals are discussed in \cite
{PeccatiTaqqu2011} and \cite{Surgailis1984}. Before stating a result
for the expected value of such a product, we introduce some notation.
The function $\bigotimes_{l=1}^m f_l\dvtx X^{\sum n_l}\rightarrow\R$ is given by
\[
\Biggl(\bigotimes_{l=1}^m
f_l\Biggr) \bigl(z_1^{(1)},
\ldots,z_{n_1}^{(1)},\ldots, z_1^{(m)},
\ldots, z_{n_m}^{(m)}\bigr)=\prod_{l=1}^m
f_l\bigl(z_1^{(l)},\ldots,z_{n_l}^{(l)}
\bigr).
\]

\begin{definition}
Let $\Pi_{n_1,\ldots,n_m}$ be the set of all partitions of the set of
variables $z_1^{(1)},\ldots, z_{n_m}^{(m)}$ such that two variables
$z^{(l)}_{i}$ and $z^{(l)}_{j}$ with $i\neq j$ but the same upper index
$(l)$ are always in different blocks, and such that every block
includes at least two variables.
\end{definition}
In this definition, we think of variables as combinatorial objects and
partition a set of them. This is slightly different from the approach
in \cite{PeccatiTaqqu2011}, where the variables are numbered, and the
partitions are defined for a set of numbers. Observe that by definition
each block of $\pi\in\Pi_{n_1,\ldots,n_m}$ has at least two and at most
$m$ variables each of them with different upper index $(l)$.
Subsequently also the following subset of $\Pi_{n_1,\ldots,n_m}$ will
play a central role.
%
\begin{definition}
Let $\overline{\Pi}_{n_1,\ldots,n_m}$ be the set of all partitions $\pi
\in\Pi_{n_1,\ldots,n_m}$ such that for any decomposition of $\{ 1,
\ldots, m\}$ into two disjoint nonempty sets $M_1, M_2$ there are $l_1
\in M_1, l_2 \in M_2$ and two variables $z_i^{(l_1)}, z_j^{(l_2)}$
which are in the same block of $\pi$.
\end{definition}

By $|\pi|$ we denote the number of blocks of the partition $\pi$. For
every partition $\pi\in\Pi_{n_1,\ldots,n_m}$ we define the function
$(\bigotimes_{l=1}^m f_l)_\pi\dvtx  X^{|\pi|}\rightarrow\R$ by replacing all
variables of $\bigotimes_{l=1}^m f_l$ that belong to the same block of $\pi
$ by a new common variable. The order of the new variables does not
matter since we always integrate over all variables.

Let us recall that ${\mathcal S}(X^k)$ stands for the set of simple
functions. These are all $f\in L^2_s(X^k)$ that are zero on all
diagonals, are constant on a finite number of Cartesian products, and
vanish everywhere else. For the product of multiple Wiener--It\^o
integrals of such functions the following proposition holds; see
Corollary 7.2 in~\cite{PeccatiTaqqu2011}.
%
\begin{proposition}\label{prop:ProductFormula}
Let $f_l\in\mathcal{S}(X^{n_l})$ for $l=1,\ldots,m$. Then
%
%
\begin{equation}
\label{eqn:ProdSimple} \qquad\E\prod_{l=1}^m
I_{n_l}(f_l)= \sum_{\pi\in\Pi_{n_1,\ldots,n_m}} \int_{X^{|\pi|}}\Biggl(\bigotimes_{l=1}^m
f_l\Biggr)_\pi(y_1,\ldots,y_{|\pi|})
\,d\mu(y_1,\ldots,y_{|\pi|}).
\end{equation}
\end{proposition}
As a consequence of Proposition 3.1 in \cite{Surgailis1984}, equation
(\ref{eqn:ProdSimple}) is also true for $f_l\in L_s^2(X^k)$, $l=1,\ldots,m$, satisfying
%
%
\begin{equation}
\label{eqn:ConditionProduct} \Biggl(\bigotimes_{l=1}^mf_l
\Biggr)_\pi\in L^2\bigl(X^{|\pi|}\bigr)
\end{equation}
for all partitions $\pi$ of the set of variables such that all
variables of a function are in different blocks. For some classes of
functions $f_l$ it is obvious that (\ref{eqn:ConditionProduct}) holds,
for example, if the $f_l$ are bounded and have a support of finite
measure. But in general it is difficult to verify condition (\ref
{eqn:ConditionProduct}).

In order to avoid this problem, we approximate a general $U$-statistic by
a sequence of $U$-statistics, whose kernels are simple functions, and
apply Proposition~\ref{prop:ProductFormula}. Afterward, we extend our
results to the original $U$-statistic. From now on, we assume that $F$ is
an absolutely convergent $U$-statistic.

Because of $f\in L_s^1(X^k)$, there exists a sequence $
(f^{(n)} )_{n\in\mathbb{N}}$ of functions in $\mathcal{S}(X^k)$
such that $|f^{(n)}|\leq|f|$ $\mu^{k}$-almost everywhere and $
(f^{(n)} )_{n\in\mathbb{N}}$ converges to $f$ $\mu^k$-almost
everywhere on $X^k$. We define $U$-statistics $F^{(n)}, n\in\N$, by
\[
F^{(n)}=\sum_{(x_1,\ldots,x_k)\in\eta_{\neq}^k}f^{(n)}(x_1,
\ldots,x_k).
\]
Since $ (f^{(n)} )_{n\in\mathbb{N}}$ converges $\mu^k$-almost
everywhere on $X^k$ to $f$,
\[
\lim_{n\rightarrow\infty}f^{(n)}(x_1,
\ldots,x_k)=f(x_1,\ldots,x_k)\qquad\mbox{for all
}(x_1,\ldots,x_k)\in\eta_{\neq}^k
\]
holds with probability $1$. Furthermore, the absolute convergence of
$F$ implies
\[
\bigl|F^{(n)}\bigr|\leq\sum_{(x_1,\ldots,x_k)\in\eta^{k}_{\neq
}}\bigl|f^{(n)}(x_1,
\ldots,x_k)\bigr|\leq\sum_{(x_1,\ldots,x_k)\in\eta^{k}_{\neq
}}\bigl|f(x_1,
\ldots,x_k)\bigr|\in L^2(\P).
\]
Hence, $ (F^{(n)} )_{n\in\mathbb{N}}$ converges almost surely
to $F$, and the dominated convergence theorem implies even convergence
in $L^1(\P)$ and $L^2(\P)$. Moreover, $F^{(n)}\in L^2(\P)$, and every
$F^{(n)}$ has a Wiener--It\^o chaos expansion with kernels $f_i^{(n)}$
that are simple functions since integration over a variable of a simple
function leads to a simple function.

The fact that the kernels of $F^{(n)}$ are simple functions brings us
in the position to use Proposition \ref{prop:ProductFormula} to
evaluate $R_{ij}$ and $\tilde{R}_{i}$ for $i,j=1,\ldots,k$.
We start by estimating $R_{ii}$. By (\ref{eqn:ProdSimple}), we have
\begin{eqnarray*}
R_{ii} &=& \int_{X^2} \E I_{i-1}
\bigl(f^{(n)}_i(s,\cdot)\bigr)^2
I_{i-1}\bigl(f^{(n)}_i(t,\cdot)
\bigr)^2 \,d\mu(s,t) - \bigl[ (i-1)!\bigl\|f^{(n)}_i
\bigr\|^2 \bigr]^2
\\
&= & \sum_{\pi\in\Pi_{i-1,i-1,i-1,i-1}} \int_{X^{|\pi|+2}}
\bigl(f^{(n)}_i(s,\cdot)\otimes f^{(n)}_i(s,
\cdot)\\
&&\hspace*{104pt}{}\otimes f^{(n)}_i (t,\cdot)\otimes
f^{(n)}_i (t,\cdot) \bigr)_\pi(y_1,\ldots,y_{|\pi|}) \\
&&\hspace*{104pt}d\mu(y_1,
\ldots,y_{|\pi|}, s,t)
\\
&&{} - \bigl[ (i-1)!\bigl\|f^{(n)}_i\bigr\|^2
\bigr]^2.
\end{eqnarray*}
The sum over those partitions of $\Pi_{i-1,i-1,i-1,i-1}$ such that
every block contains only variables of the first pair or of the second
pair of functions leads exactly to $[ (i-1)!\|f^{(n)}_i\|^2 ]^2$. These
partitions cancel out with the minus term and we denote the remaining
partitions by $\tilde{\Pi}_{i-1,i-1,i-1,i-1}$. Hence,
\begin{eqnarray*}
 R_{ii}& =& \sum_{\tilde{\pi}\in\tilde{\Pi
}_{i-1,i-1,i-1,i-1}} \int
_{X^{|\tilde{\pi}|+2}} \bigl(f^{(n)}_i(s,\cdot)\otimes
f^{(n)}_i(s,\cdot)\\
&&\hspace*{104pt}{}\otimes f^{(n)}_i(t,
\cdot) \otimes f^{(n)}_i(t,\cdot) \bigr)_{\tilde{\pi}}
 (y_1,\ldots,y_{|\tilde{\pi}|})
\\
&&\hspace*{103pt}d
\mu(y_1,\ldots,y_{|\tilde{\pi}|},s,t)
\\
 &\leq& \sum_{\tilde{\pi}\in\tilde{\Pi}_{i-1,i-1,i-1,i-1}} \int_{X^{|\tilde{\pi}|+2}}\bigl |
\bigl(f^{(n)}_i(s,\cdot) \otimes f^{(n)}_i(s,
\cdot)  \\
&&\hspace*{104pt}{}\otimes f^{(n)}_i(t,\cdot)\otimes
f^{(n)}_i(t,\cdot )\bigr)_{\tilde{\pi}} (y_1,\ldots,y_{|\tilde{\pi}|})\bigr| \\
&&\hspace*{103pt}{}d
\mu(y_1,\ldots,y_{|\tilde{\pi}|},s,t).
\end{eqnarray*}
In order to simplify our notation, we include $s$ and $t$ into the
partitions by adding two blocks generating $s$ and $t$ to the old
partition $\tilde{\pi}$ and obtain a new partition $\pi\in\Pi
_{i,i,i,i}$. By definition of $\tilde{\pi}$, $\pi$ has at least one
block including variables $z^{(l_1)}_{i_1}$ and $ z^{(l_2)}_{i_2}$,
$l_1\in \{1,2 \}$, $l_2\in \{3,4 \}$. By construction
of $\pi$, there are also blocks including variables of the first two
functions and of the last two functions. Altogether, this implies $\pi
\in\overline{\Pi}_{i,i,i,i}$. Since each $\tilde{\pi}\in\tilde{\Pi
}_{i-1,i-1,i-1,i-1}$ leads to a different $\pi\in\overline{\Pi
}_{i,i,i,i}$, we obtain the upper bound
\begin{eqnarray*}
&&R_{ii} \leq \sum_{\pi\in\overline{\Pi}_{i,i,i,i}} \int
_{X^{|\pi|}} \bigl|\bigl(f^{(n)}_i \otimes
f^{(n)}_i \otimes f^{(n)}_i \otimes
f^{(n)}_i\bigr)_{\pi
}(y_1,
\ldots,y_{|\pi|})\bigr| \\
&&\hspace*{86pt}d\mu(y_1,\ldots,y_{|\pi|}).
\end{eqnarray*}
In the very same way, we obtain an upper bound for $R_{ij},i\neq j$. By
(\ref{eqn:ProdSimple}), it follows
\begin{eqnarray*}
 R_{ij}
& = & \int_{X^2} \E \bigl[ I_{i-1}
\bigl(f^{(n)}_i(s,\cdot)\bigr) I_{j-1}
\bigl(f^{(n)}_j(s,\cdot)\bigr) I_{i-1}
\bigl(f^{(n)}_i(t,\cdot )\bigr)\\
&&\hspace*{150pt}{}\times I_{j-1}\bigl(f^{(n)}_j(t,\cdot)\bigr) \bigr] \,d\mu(s,t)
\\
& = & \sum_{\tilde{\pi}\in\Pi_{i-1,j-1,i-1,j-1}} \int_{X^{|\tilde
{\pi}|+2}}
\bigl(f^{(n)}_{i}(s,\cdot) \otimes f^{(n)}_{j}(s,
\cdot) \\
&&\hspace*{104pt}{}\otimes f^{(n)}_{i}(t,\cdot) \otimes
f^{(n)}_{j}(t,\cdot)\bigr)_{\tilde{\pi}}(y_1,\ldots,y_{|\tilde{\pi}|})\\
&&\hspace*{103pt}{}d\mu(y_1,
\ldots,y_{|\tilde{\pi
}|},s, t)
\\
& \leq&  \sum_{\pi\in\overline{\Pi
}_{i,j,i,j}} \int
_{X^{|\pi|}}\bigl |\bigl(f^{(n)}_{i} \otimes
f^{(n)}_{j} \otimes f^{(n)}_{i} \otimes
f^{(n)}_{j}\bigr)_\pi(y_1,\ldots,y_{|\pi|})\bigr| \,d\mu(y_1,\ldots,y_{|\pi|}).
\end{eqnarray*}
Since $i\neq j$, there exist no $\tilde{\pi}\in\Pi_{i-1,j-1,i-1,j-1}$
with blocks including either variables of the first two or last two
functions. Hence, we obtain partitions $\pi\in\overline{\Pi}_{i,j,i,j}$
by the same construction as for $R_{ii}$ and obtain an upper bound by
summing over $\overline{\Pi}_{i,j,i,j}$.

The last step is to estimate $\tilde R_i$. Here, we have in a similar way
\begin{eqnarray*}
 \tilde R_{i} &=& \int_{X} \E I_{i-1}
\bigl(f^{(n)}_i(s,\cdot)\bigr)^4 \,d\mu(s)
\\
& =& \sum_{\tilde{\pi}\in\Pi_{i-1,i-1,i-1,i-1}} \int_{X^{|\tilde
{\pi}|+1}}
\bigl(f_i^{(n)}(s,\cdot) \otimes f_i^{(n)}(s,
\cdot)\\
&&\hspace*{103pt}{} \otimes f_i^{(n)}(s,\cdot) \otimes
f_i^{(n)}(s,\cdot) \bigr)_{\tilde{\pi}}(y_1,\ldots,y_{|\tilde{\pi}|}) \\
&&\hspace*{101pt}{}d\mu(y_1,
\ldots,y_{|\tilde{\pi}|},s)
\\
& \leq& \sum_{\pi\in\overline{\Pi}_{i,i,i,i}} \int_{X^{|\pi|}} \bigl|
\bigl(f_i^{(n)} \otimes f_i^{(n)}
\otimes f_i^{(n)} \otimes f_i^{(n)}
\bigr)_{\pi
}(y_1,\ldots,y_{|\pi|})\bigr| \,d
\mu(y_1,\ldots,y_{|\pi|}).
\end{eqnarray*}
In this case, it is immediate that we obtain a partition $\pi\in
\overline{\Pi}_{i,i,i,i}$ by adding $s$ to a partition $\tilde{\pi}\in
\Pi_{i-1,i-1,i-1,i-1}$.
Thus $R_{ij}$ and $\tilde R_i$ are bounded by the same expressions.

Now it remains to estimate the kernels $f_i^{(n)}$. From Lemma \ref
{Lemma:Kernels}, it follows that
\begin{eqnarray*}
&&\bigl|f_i^{(n)} (y_1,\ldots,y_i)\bigr| \\
&&\qquad\leq
\pmatrix{k
\cr
i} \int_{X^{k-i}}\bigl|f^{(n)}
(y_1,\ldots,y_i,x_1,\ldots,x_{k-i})\bigr|
\,d\mu (x_1,\ldots,x_{k-i}).
\end{eqnarray*}
We obtain the following expression as an upper bound for $R_{ij}$ and
$\tilde R_{i}$. With $M_{ij} (\cdot) $ defined by
%
%
\begin{eqnarray}\label{def:Mij}
M_{ij}(g) &=& \pmatrix{k
\cr
i}^2\pmatrix {k
\cr
j}^2\nonumber\\
&&{}\times \sum_{\pi\in\overline{\Pi
}_{i,j,i,j}} \int
_{X^{|\pi|+4k-2i-2j}}\bigl |\bigl( g\bigl(\cdot, x^{(1)}_1,
\ldots,x^{(1)}_{k-i}\bigr) \nonumber\\
&&\hspace*{110pt}{}\otimes
 g\bigl(\cdot, x^{(2)}_1,
\ldots,x^{(2)}_{k-j}\bigr)
\nonumber
\\[-8pt]
\\[-8pt]
\nonumber
&&\hspace*{110pt}{}\otimes g\bigl(\cdot,
x^{(3)}_1,\ldots,x^{(3)}_{k-i}\bigr)\\
&&\hspace*{110pt}{}\otimes g\bigl(\cdot, x^{(4)}_1,\ldots,x^{(4)}_{k-j}
\bigr) \bigr)_{\pi} (y_1,\ldots,y_{|\pi|})\bigr| \nonumber\\
&&\hspace*{106pt}{} d\mu
\bigl(x_1^{(1)},\ldots,x_{k-j}^{(4)},
y_1, \ldots, y_{|\pi|}\bigr),\nonumber
\end{eqnarray}
where $\pi$ acts on the first $i$, respectively, $j$ variables of $g\dvtx X^k \to\R$,
we have
\[
R_{ij}\leq M_{ij}\bigl(f^{(n)}\bigr)  \quad\mbox{and}\quad
 \tilde{R}_i\leq M_{ii}\bigl(f^{(n)}\bigr)  \qquad\mbox{for } 1\leq i,j\leq k.
\]
Since in the definition of $M_{ij}$ every block of a partition $\pi\in
\overline{\Pi}_{i,j,i,j}$ has at least two elements, the integration in
(\ref{def:Mij}) runs over at most $4k-i-j$ variables. For $i=j=1$ the
only partition in $\overline{\Pi}_{i,j,i,j}$ is the partition with one
block and the integration runs over $4k-3$ variables. This observation
will be important in Section~\ref{sec:GeometricUstatistic}.

Combining our bounds for $R_{ij}$ and $\tilde{R}_i$ with Theorem \ref
{thm:GeneralCLT} yields:
%
\begin{lemma}
Suppose $F^{(n)}=\sum_{\eta^k_{\neq}} f^{(n)}$ is a $U$-statistic of
order $k$ with $f^{(n)}\in\mathcal{S}(X^k)$ and $N$ is a standard
Gaussian random variable. Then
\[
d_W \biggl(\frac{F^{(n)}-\E F^{(n)}}{\sqrt{\V F^{(n)}}},N \biggr)\leq 2k^{7 /2} \sum
_{1\leq i \leq j\leq k}\frac{\sqrt {M_{ij}(f^{(n)})}}{\V F^{(n)}}.
\]
\end{lemma}
Together with the fact that $M_{ij}(f^{(n)})\leq M_{ij}(f)$ since $|
f^{(n)}| \leq| f|$ and the triangle inequality for the Wasserstein
distance, we obtain
\begin{eqnarray*}
&& d_W \biggl(\frac{F-\E F}{\sqrt{\V F}},N \biggr)
\\
& &\qquad\leq d_W \biggl(\frac{F-\E F}{\sqrt{\V F}},\frac{F^{(n)}-\E
F^{(n)}}{\sqrt{\V F^{(n)}}}
\biggr)+d_W \biggl(\frac{F^{(n)}-\E
F^{(n)}}{\sqrt{\V F^{(n)}}},N \biggr)
\\
& &\qquad \leq d_W \biggl(\frac{F-\E
F}{\sqrt{\V F}},\frac{F^{(n)}-\E F^{(n)}}{\sqrt{\V F^{(n)}}} \biggr)+
2k^{7 /2} \sum_{1\leq i \leq j\leq k}\frac{\sqrt{M_{ij}(f)}}{\V F^{(n)}}.
\end{eqnarray*}
By the definition of the Wasserstein distance and some straightforward
computations, it follows that
\begin{eqnarray*}
 d_W \biggl(\frac{F-\E F}{\sqrt{\V F}},\frac{F^{(n)}-\E F^{(n)}}{\sqrt {\V F^{(n)}}} \biggr)& \leq&\E
\biggl\llvert \frac{F-\E F}{\sqrt{\V F}}-\frac
{F^{(n)}-\E F^{(n)}}{\sqrt{\V F^{(n)}}}\biggr\rrvert
\\
& \leq& \frac{\E
|F^{(n)}-F+\E F - \E F^{(n)} |}{\sqrt{\V F^{(n)}}} \\
&&{}+\biggl\llvert \frac{\E|F-\E
F|}{\sqrt{\V F}}-\frac{\E|F-\E F|}{\sqrt{\V F^{(n)}}}
\biggr\rrvert.
\end{eqnarray*}
Because of the convergence of $(F^{(n)})_{n\in\N}$ to $F$ in $L^1(\P)$
and $L^2(\P)$, the right-hand side vanishes for $n\rightarrow\infty$,
and we get our main result.

\begin{theorem}\label{thm:CLT}
Suppose $F $ is an absolutely convergent $U$-statistic of order $k$, and
$N$ is a standard Gaussian random variable. Then
\[
d_W \biggl(\frac{F-\E F}{\sqrt{\V F}},N \biggr)\leq 2k^{7 /2} \sum
_{1\leq i \leq j\leq k}\frac{\sqrt{M_{ij}(f)}}{\V F}
\]
with $M_{ij}(f)$ defined in (\ref{def:Mij}).
\end{theorem}


\section{Geometric $U$-statistics}\label{sec:GeometricUstatistic}


\subsection{Central limit theorems for geometric $U$-statistics}
In this section, we assume that our intensity measure has the form $\mu
(\cdot)=\lambda\theta(\cdot)$ with a $\sigma$-finite nonatomic measure
$\theta(\cdot)$ and $\lambda\geq1$. We are interested in the behavior
of the $U$-statistic $F$ for $\lambda\rightarrow\infty$.

\begin{definition}
A $U$-statistic $F = \sum_{\eta^k_{\neq}} f$ is a \textit{geometric}
$U$-statistic if it satisfies
\[
f(x_1,\ldots,x_k)=g(\lambda)\tilde{f}(x_1,
\ldots,x_k)
\]
with $g\dvtx \R\to\R$, and with $\tilde{f}\dvtx X^k\to\R$ not depending on
$\lambda$.
\end{definition}
In the case that $g=1$ and $f=\tilde{f}$, the value of $F$ for a given
realization of the Poisson point process is only determined by the
geometry of the points and does not depend on the intensity rate
$\lambda$ of the underlying process. The term ``geometric'' is used to emphasize this behavior. We slightly generalize this
property by allowing our geometric $U$-statistics to have an intensity
related scaling factor since we always consider standardized random
variables, where the scaling factor is cancelled out.

By $\tilde{M}_{ij}$ we denote the value of $M_{ij}(\tilde f)$, which is
defined in (\ref{def:Mij}), for $\lambda=1$. With this notation, the
following central limit theorem holds:

\begin{theorem}\label{thm:CLTgeometric}
Suppose $F $ is an absolutely convergent geometric $U$-statistic of order
$k$ with $\|f_1\|>0$ and $N$ is a standard Gaussian random variable.
Then
\[
\lim_{\lambda\rightarrow\infty}\frac{\V F}{\lambda^{2k-1}
g(\lambda)^2} =\underbrace{ k^2
\int_{X} \biggl( \int_{X^{k-1}} \tilde
f(y,x_1,\ldots,x_{k-1}) \,d\theta(x_1, \ldots,
x_{k-1}) \biggr)^2 \,d\theta(y)}_{=:\tilde{V}}
\]
with $\tilde{V} >0$, and
%
%
\begin{equation}
\label{eqn:CLTgeometric} d_W \biggl(\frac{F-\E F}{\sqrt{\V F}},N \biggr) \leq
\lambda^{-1/ 2} \biggl( 2k^{7 /2}\sum
_{1\leq i\leq j\leq
k}\frac{\sqrt{\tilde{M}_{ij}}}{\tilde{V}} \biggr)
\end{equation}
for $\lambda\geq1$.
\end{theorem}
The main feature of this theorem is that the term in brackets is
independent of~$\lambda$, which means that for $\lambda\to\infty$ the
distance to the Gaussian distribution tends to zero at a rate $\lambda
^{- 1/2}$.

\begin{pf*}{Proof of Theorem \ref{thm:CLTgeometric}}
Because we are interested in the standardized variable $(F-\E F)/ \sqrt {\V F}$ which is independent of $g(\lambda)$, w.l.o.g. we put $g(\lambda
)=1$ and $\tilde f = f$. From formula (\ref{eq:GeneralVariance}) we infer
\begin{eqnarray*}
\V F& =& \sum_{i=1}^k
\lambda^{2k-i} i!\pmatrix{k
\cr
i}^2
\\
&&\hspace*{14pt}\times{} \int_{X^i} \biggl( \int
_{X^{k-i}} f(y_1,\ldots,y_i,x_1,
\ldots,x_{k-i}) \,d\theta(x_1, \ldots, x_{k-i})
\biggr)^2 \\
&&\hspace*{42pt}{}d\theta(y_1, \ldots, y_{i}),
\end{eqnarray*}
which means that the variance is a polynomial of degree $2k-1$ with the
leading term $\tilde V \lambda^{2k-1} = \| f_1\|^2 >0 $ and
$\V F\geq\tilde{V} \lambda^{2k-1}$.

As previously mentioned, the integration in $M_{ij}(f)$ runs over at
most $4k-i-j\leq4k-3$ variables for $(i,j)\neq(1,1)$ and $4k-3$
variables for $(i,j)=(1,1)$, and we see that
\[
M_{ij}(f)\leq\tilde{M}_{ij} \lambda^{4k-3}
\]
for $\lambda\geq1$. Hence, Theorem \ref{thm:CLT} leads directly to
(\ref{eqn:CLTgeometric}).
\end{pf*}

The assumption $\| f_1\| >0$ cannot be easily dispensed as can be seen
from the following example:

\begin{example*}
Let $\eta$ be a Poisson process on $[-1,1]$ with intensity measure the
Lebesgue measure times intensity\vadjust{\goodbreak} $\lambda>0$. We define the $U$-statistic
$F=\sum_{(x_1,x_2)\in\eta^2_{\neq}}f(x_1,x_2)$ with
\[
f(x_1,x_2)= %
\cases{1, &\quad $x_1x_2
\geq0,$\vspace*{2pt}
\cr
-1, &\quad $x_1x_2<0$.} %
\]
Obviously, we obtain $f_1(y)=0$. It is possible to rewrite $F$ as
$F=L(L-1)+R(R-1)-2LR$
where $L$ and $R$ are the number of points in $[-1, 0]$ and $[0,1]$,
respectively. This brings us in the position to compute the moments.
Elementary calculations show that the variance equals $8 \lambda^2$,
and the third moment of $F$ is $64\lambda^3$. Thus the third moment of
$(F - \E F)/ \sqrt{\V F}$ tends to a constant and hence is too large
for convergence of $F$ to a Gaussian distribution.
By a technical computation of all moments, using the product formula
for multiple Wiener--It\^o integrals, for example, and the method of
moments, it can be shown that $\sqrt{2}(F-\E F)/\sqrt{\V F}$ follows a
centered chi-square distribution with one degree of freedom as $\lambda
\rightarrow\infty$.
\end{example*}

In the special case $\mu(X)=\lambda\theta(X) < \infty$, it is possible
to approximate the Poisson point process $\eta$ by a binomial point
process, that consists of a fixed number of independently distributed
points with the probability measure $\theta(\cdot)/\theta(X)$. If we
sum over $k$-tuples of distinct points of the binomial point process
instead of a Poisson point process, we obtain a classical $U$-statistic.
This well-known class of random variables satisfies a similar central
limit theorem as above with a rate of convergence; see \cite
{Hoeffding1948,GramsSerf,KoroljukBorovskich,Lee1990}. Although both
results are similar, it seems to be difficult to prove one result by
the other, especially with keeping rates of convergence.

For classical $U$-statistics the so-called Hoeffding decomposition which
is closely related to the Wiener--It\^o chaos expansion plays a crucial
role. In the recent paper by Lachi\'eze-Rey and Peccati \cite
{LachPec2}, this decomposition is applied to $U$-statistics of Poisson
point processes which yields a representation similar to the
Wiener--It\^o chaos expansion. Combining this with the result of Dynkin
and Mandelbaum \cite{DynkinMandelbaum1983}, the authors derive our
Theorem \ref{thm:CLTgeometric} for the case $\mu(X)<\infty$ (without
rates of convergence).
They also prove noncentral limit theorems for the case that some of
the first kernels of the chaos expansion of a $U$-statistic vanish, which
allows one to deal with situations as in the previous example.

In Sections~\ref{subsec:PHP} and \ref{subsec:ConvexHull}, we apply
Theorem \ref{thm:CLTgeometric} to problems from stochastic geometry.
In the recent paper \cite{Decreusefondetal2011} the underlying result
from \cite{Peccatietal2010} is used to derive a central limit theorem
for the number random simplices on a torus. This problem exactly fits
in the framework of geometric $U$-statistics, and some of the results can
also be obtained by using Theorem \ref{thm:CLTgeometric}.


\subsection{Central limit theorems for Poisson hyperplanes}\label{subsec:PHP}
We use Theorem~\ref{thm:CLTgeometric} to establish central limit
theorems for Poisson hyperplane processes. Let $\eta$ be a
Poisson
process on the space $\mathcal{H}$ of all hyperplanes in $\R^d$ with an
intensity measure of the form $\mu(\cdot)=\lambda\theta(\cdot)$ with
$\lambda\in\R^+$ and a $\sigma$-finite nonatomic\vadjust{\goodbreak} measure $\theta$. The
Poisson hyperplane process is only observed in a compact convex window
$W\subset\R^d$ with interior points. Thus, we can view $\eta$ as a
Poisson process on the set $[W]$ defined by
\[
[W]= \{h\in\mathcal{H}\dvtx h\cap W\neq\varnothing \}.
\]

Given the hyperplane process $\eta$, we investigate the $(d-k)$-flats
in $W$ which occur as the intersection of $k$ hyperplanes of $\eta$. In
particular, we are interested in the sum of their $i$th intrinsic
volumes given by
\[
\label{eqn:Intersections} \Phi^k_i(W)=
\frac{1}{k!} \sum_{(h_1,\ldots,h_k)\in\eta^{k}_{\neq}} V_i(
h_1\cap\cdots\cap h_k\cap W )
\]
for $i=0, \ldots, d-k$ and $k=1,\ldots,d$.
For the definition of the $i$th intrinsic volume $V_i(\cdot)$ we refer
to \cite{Schn4}. We remark that $V_0(K)$ is the Euler characteristic of
the set~$K$, and that $V_{n} (K) $ of an $n$-dimensional convex set~$K$
is the Lebesgue measure $\Lambda_n (K)$. Thus $\Phi^k_0$ is the number
of $(d-k)$-flats in $W$, and $\Phi^k_{d-k}$ is their $(d-k)$-volume.
To ensure that the expectations of these random variables are neither 0
nor infinite, we assume that:
\begin{itemize}
\item$0<\theta([W])<\infty$;
\item$2\leq k \leq d$ independent random hyperplanes on $[W]$ with
probability measure $\theta(\cdot)/\theta([W])$ intersect in a
$(d-k)$-flat almost surely and their intersection flat hits the
interior of $W$ with positive probability.
\end{itemize}
For example, these conditions are satisfied if the hyperplane process
is stationary and the directional distribution is not concentrated on a
great subsphere.

The fact that the summands in the definition of $\Phi_i^k$ are bounded
and have a bounded support makes sure that the fourth moments in
$M_{ij}(\cdot)$ are finite, and we can apply Theorem \ref{thm:CLTgeometric}:
%
\begin{theorem}
Let $N$ be a standard Gaussian random variable. Then constants $c_{\Phi
}(W,k,i)$ exist such that
\[
d_W \biggl(\frac{\Phi_i^k(W)-\E  \Phi_i^k(W)}{\sqrt{\V\Phi_i^k
(W)}},N \biggr)\leq c_{\Phi}(W,k,i)
\lambda^{-1/ 2}
\]
for $\lambda\geq1$, $i=0, \ldots, d-k$ and $k=1,\ldots,d$.
\end{theorem}

Furthermore, the asymptotic variances are given by
\begin{eqnarray*}
 &&\lim_{\lambda\to\infty} \frac{\V\Phi_i^k(W)
}{\lambda^{2k-1}}
\\
&&\qquad= \frac{1}{(k-1)!^2}\\
&&\qquad\quad{}\times \int_{[W]} \biggl(
\int_{[W]^{k-1}} V_i( h \cap h_1\cap\cdots
\cap h_{k-1}\cap W ) \,d\theta(h_1, \ldots, h_{k-1})
\biggr)^2 \,d\theta(h).
\end{eqnarray*}

Similar results have first been derived by Paroux \cite{Paroux1998},
and by Heinrich \cite{Heinrich2009} and Heinrich, Schmidt and Schmidt
\cite{HeinrichSchmidtSchmidt2006} using Hoeffding's decomposition of
classical $U$-statistics. Schulte and Th\"ale \cite{SchulteThaele2010}
used the Wiener--It\^o chaos expansion to compute the moments and
cumulants and to formulate central limit theorems for the surface area
of Poisson hyperplanes in an increasing window. In their recent
paper~\cite{SchulteThaele2012} this approach is further refined to obtain
point process convergence for the intrinsic volumes of the intersection
process of Poisson $k$-flats in the unit ball.


\subsection{Convex hulls of random points}\label{subsec:ConvexHull}

In the following, we assume that the Poisson point process $\eta$ has
an intensity-measure of the form $\mu(\cdot)=\lambda\Lambda_d(\cdot
\cap K)$, $\lambda\geq1$, where $\Lambda_d$ is Lebesgue measure, and
$K\subset\R^d$ a compact convex set with $\Lambda_d(K)=1$. If we
integrate with respect to $\Lambda_d$, we omit the measure in our notation.

We consider the following functional related to Sylvester's problem:
\[
H=\sum_{(x_1,\ldots,x_k)\in\eta^k_{\neq}}h(x_1,\ldots,x_k)
\]
with
\[
h(x_1,\ldots,x_k)=\1 \bigl(x_1,
\ldots,x_k \mbox{ are vertices of }\operatorname{conv}(x_1,
\ldots,x_k)\bigr),
\]
which counts the number of $k$-tuples of the process such that every
point is a vertex of the convex hull, that is, the number of $k$-tuples
in convex position. The expected value of $H$ is then given by
\begin{eqnarray*}
\E H &=& \lambda^k \P\bigl(X_1,\ldots,X_k
\mbox{ are vertices of }\operatorname{conv}(X_1,\ldots,X_k)\bigr) =
\lambda^k p^{(k)}(K),
\end{eqnarray*}
where $X_1,\ldots,X_k$ are independent random points chosen according
to the uniform distribution on $K$.

The question to determine the probability $p^{(k)}(K)$ that $k$ random
points in a convex set $K$ are in convex position has a long history;
see, for example, the more recent developments by B{\'a}r{\'a}ny \cite
{Bar5,Bar6} and Buchta \cite{Bu12}.
In our setting, the function $H$ is an estimator for the probability
$p^{(k)}(K)$, and we are interested in distributional properties of
this estimator.

The asymptotic behavior of $\V H$ is determined by
\[
\tilde H = \lim_{\lambda\rightarrow\infty}\frac{\V H}{\lambda^{2k-1} } = k^2
\int_K \biggl(\int_{K^{k-1}}h(y,x_1,
\ldots,x_{k-1}) \,dx_1 \cdots dx_{k-1}
\biggr)^2 \,dy.
\]
By the Cauchy--Schwarz inequality, because $\Lambda_d(K)=1$ and $h^2 =
h$, we obtain
\begin{eqnarray*}
\tilde H &\leq& k^2 \int_K \int
_{K^{k-1}}h(y,x_1,\ldots,x_{k-1})^2
\,dx_1 \cdots dx_{k-1} \,dy
\\
&=& k^2 p^{(k)}(K)
\end{eqnarray*}
and
\begin{eqnarray*}
k^2 p^{(k)}(K)^2 &=& k^2 \biggl(
\int_{K^k}h(x_1,\ldots,x_k) \,d
x_1 \cdots dx_{k} \biggr)^2
\\
&\leq& k^2 \int_K \biggl( \int
_{K^{k-1}} h(x_1,x_2,\ldots,x_k)
\,dx_2 \cdots dx_k \biggr)^2 \,dx_1 =
\tilde H.
\end{eqnarray*}
Together with Theorem \ref{thm:CLTgeometric}, we immediately get the
following result showing that the estimator $H$ is asymptotically Gaussian:

\begin{theorem}
Let $N$ be a standard Gaussian random variable. Then there exists a
constant $C$ such that
\[
d_W \biggl(\frac{H-\E  H}{\sqrt{\V H}},N \biggr)\leq C \lambda^{-1/
2}.
\]
Furthermore
$ \V H = \lambda^{2k-1} \tilde H (1+ O(\lambda^{-1}))$ as $\lambda\to
\infty$ with
\[
k^2 p^{(k)}(K)^2 \leq\tilde H \leq
k^2 p^{(k)}(K).
\]
\end{theorem}


\section{Local $U$-statistics}\label{sec:LocalUstatistic}


\subsection{Central limit theorems for local $U$-statistics}
For a geometric $U$-statistic the function $f$ is [up to the scaling
factor $g(\lambda)$] independent of $\lambda$. Now we allow that $f$ is
influenced by $\lambda$ in a more intricate way, but we assume that a
$k$-tupel of points is only in the support of $f$ if the points are
close together.

From now on, let $X$ be a metric space, and denote by $B(y, r)$ the
ball with center $y$ and radius $r$. Again, we assume that the
intensity measure $\mu$ has the form $\mu(\cdot)=\lambda\theta(\cdot)$
with $\lambda\geq1$ and a $\sigma$-finite nonatomic measure $\theta
(\cdot)$ on $X$. We denote the diameter of $A\subset X$ by $\operatorname{diam}(A)$.

\begin{definition}
A $U$-statistic $F = \sum_{\eta^k_{\neq}} f$ is a \textit{local}
$U$-statistic if it satisfies
%
%
\begin{equation}
\label{eqn:flocal1} f(x_1,\ldots,x_k)=0\qquad  \mbox{if }  \operatorname{diam}\bigl( \{x_1,\ldots,x_k \}\bigr) > \delta.
\end{equation}
\end{definition}

Note that in general $\delta$ may depend on $\lambda$.
We denote the $L^2$-norm on $X^i$ with respect to the measure $\theta
(\cdot)$ by $\|\cdot\|_{\theta}$. Now we can rephrase Theorem \ref
{thm:CLT} for local $U$-statistics as follows:\vadjust{\goodbreak}

\begin{theorem}\label{thm:CLTlocal2}
Suppose $F $ is an absolutely convergent local $U$-statistic of order $k$
with $\|f_1\| >0$, and $N$ is a standard Gaussian random variable.
Putting $\tilde V = \| f_1 \|^2/\lambda^{2k-1}$ and $b(\delta) = \max_{y \in X} \mu(B(y, 4 \delta)) < \infty$, we have
\[
d_W \biggl(\frac{F -\E F}{\sqrt{\V F}},N \biggr)\leq c_k
\lambda^{- {3k}/2 + 1} \max\bigl\{ 1, b(\delta)^{{(3k-3)}/{2}} \bigr\}
\frac{ \| f^2\|_{\theta} } {
\tilde V}
\]
with a constant $c_k\in\R$ only depending on $k$.
\end{theorem}

\begin{pf}
Formula (\ref{eq:GeneralVariance}) yields $\V F\geq\|f_1\|^2 = \lambda
^{2k-1}\tilde V$. The estimate for $M_{ij}$ runs as follows.
Since $\pi\in\overline{\Pi}_{i,j,i,j}$, and condition (\ref
{eqn:flocal1}) forces all arguments of $f$ to be close, we can rewrite
$M_{ij}(f)$ as
\begin{eqnarray*}
M_{ij}(f)  &= & \pmatrix{k
\cr
i}^2\pmatrix{k
\cr
j}^2 \\
&&\times{}\sum_{\pi\in\overline{\Pi}_{i,j,i,j}} \int
_{X^{|\pi|+4k-2i-2j}} \bigl|\bigl( f\bigl(\cdot,x_1^{(1)},
\ldots,x_{k-i}^{(1)}\bigr)\\
&&\hspace*{107pt}{} \otimes f\bigl(\cdot,x_1^{(2)},\ldots,x_{k-j}^{(2)}
\bigr) \\
&&\hspace*{107pt}{}\otimes f\bigl(\cdot,x_1^{(3)},\ldots,x_{k-i}^{(3)}
\bigr) \\
&&\hspace*{107pt}{}\otimes f\bigl(\cdot,x_1^{(4)},\ldots,x_{k-j}^{(4)}
\bigr)\bigr)_\pi(y_1,\ldots,y_{|\pi|})\bigr |\\
&&\hspace*{101pt}{}\times \1 \bigl(\operatorname{diam} \bigl( \bigl\{ x_1^{(1)},
\ldots,x^{(4)}_{k-j}, y_1, \ldots, y_{|\pi
|}
\bigr\} \bigr) \leq4 \delta\bigr)\\
&&\hspace*{106pt}{} d\mu\bigl(x_1^{(1)}, \ldots,x^{(4)}_{k-j},
y_1, \ldots, y_{|\pi|}\bigr).
\end{eqnarray*}
By H\"older's inequality, we obtain
\begin{eqnarray*}
M_{ij} (f) &\leq& c_{ij} \sum_{\pi\in\overline{\Pi}_{i,j,i,j}}
\int_{X^{|\pi|+4k-2i-2j}} f(z_1, \ldots, z_k)
^4\\
 &&{}\times\1 \bigl(\operatorname{diam} \bigl(\{ z_1, \ldots, z_{|\pi|+4k-2i-2j} \}
\bigr) \leq4 \delta \bigr) \,d\mu(z_1, \ldots, z_{|\pi|+4k-2i-2j})
\\
&\leq& c_{ij} \bigl\| f^2\bigr\|^2 \sum
_{\pi\in\overline{\Pi}_{i,j,i,j}} b(\delta)^{|\pi| +3k-2i-2j}
\\
& = & c_{ij} \lambda^{k} \bigl\| f^2
\bigr\|_{\theta}^2  \sum_{\pi\in\overline{\Pi}_{i,j,i,j}} b(
\delta)^{|\pi| +3k-2i-2j}
\end{eqnarray*}
with a constant $c_{ij}\in\R$ depending on $i,j,k$. One should keep in
mind that $ \max(i,j) \leq| \pi| \leq i+j$ for all $\pi\in\overline
{\Pi}_{i,j,i,j}$ and that the only partition $\pi\in\overline{\Pi
}_{1,1,1,1}$ satisfies $|\pi|=1$. This leads to $|\pi|-2i-2j\leq-3$ and
\begin{eqnarray*}
2k^{7 /2} \sum_{1\leq i \leq j\leq k}\frac{\sqrt{M_{ij}(f)}}{\V F} &
\leq& c_k' \sum_{1\leq i \leq j\leq k}
\frac{ \sqrt{ \lambda^{k} \| f^2\|_{\theta}^2 \sum_{\pi\in\overline
{\Pi}_{i,j,i,j}} b(\delta)^{|\pi| +3k-2i-2j} }} {
\lambda^{2k-1} \tilde V }
\\
& \leq& c_k' \lambda^{- {3k}/2 + 1} \frac{ \| f^2\|_{\theta} } {
\tilde V }
\sum_{1\leq i \leq j\leq k} \sqrt{\sum_{\pi\in\overline{\Pi}_{i,j,i,j}}
b(\delta)^{|\pi|
+3k-2i-2j} }
\\
& \leq& c_k \lambda^{- {3k}/2 + 1} \frac{ \| f^2\|_{\theta} } {
\tilde V} \max\bigl\{ 1,
b(\delta)^{{(3k-3)}/{2}} \bigr\}
\end{eqnarray*}
with constants $c_k',c_k\in\R$ only depending on $k$. Combining this
estimate with Theorem \ref{thm:GeneralCLT} gives the claimed result.
\end{pf}
The proof rests essentially upon the fact that $F$ is a local
$U$-statistic since this allows us to rewrite $M_{ij}(f)$ such that every
function depends on all variables and to split these functions using H\"
older's inequality.


\subsection{A central limit theorem for the total edge length of a
random geometric graph}
We apply the results of the previous subsection to a problem from
random graph theory. We construct a random graph in the following way.
Let $\eta$ be a Poisson process in $X=\R^d$ with an intensity measure
of the form
\[
\mu(\cdot)=\lambda\Lambda_d(\cdot\cap W)
\]
with $\lambda\geq1$, the $d$-dimensional Lebesgue-measure $\Lambda
_d(\cdot)$ and a compact window $W \subset\R^d$ of volume $\Lambda
_d(W)=1$ containing the origin in its interior.
We regard $\eta$ as a set of points in $W$. As in (\ref{eqn:flocal1})
we connect two points $x,y\in\eta$ by an edge if
\[
\|x-y\| \leq\d=\d(\lambda).
\]
The resulting graph $G(P_\lambda,\d)$ is a random geometric graph,
sometimes called a Gilbert graph or an interval graph (for $d=1$) and a
disk graph (for $d=2$). For graph-theoretical properties of $G(P_\lambda, \d)$ we refer to \cite{PenroseRandomGeometricGraphs} and to the more
recent developments \cite{HanMak,McDi3,Muller,LachPec1}.
For our central limit theorem we take $\lambda\rightarrow\infty$ and
assume that $\d$ is small enough to ensure that
\[
\bigcap_{x \in B(0, \d)} (W+x) \supset\frac12 W.
\]

We are interested in the total edge length $L(\eta)$ of $G(P_\lambda,\d
)$ in the window $W$, which is given by
\[
L(\eta)=\frac{1}{2}\sum_{(x,y)\in\eta^2_{\neq}} g( x-y ) \1\bigl(\|
x-y\| \leq\d\bigr).
\]
Here $g\dvtx B(0,\delta) \to\R$ is some kind of measure of the length of
the edge $(x,y)$.
We assume $g \in L^2 (B(0, \delta))$ which implies that $L$ is
absolutely convergent.
The following lemma is immediate from Lemma \ref{Lemma:Kernels}.

\begin{lemma}\label{le:Leta}
$L(\eta)$ has a Wiener--It\^o chaos expansion with kernels
\[
f_1(y)=\lambda\int_{B(0, \d)} g( x) \1(y+x \in W )
\,dx,\qquad y\in W
\]
and
\[
f_2(x,y)=\tfrac{1}{2}g(x-y)\1 \bigl(\|x-y\|\leq\delta \bigr),\qquad x,y\in W.
\]
\end{lemma}

For the length of this random graph, we obtain the following central
limit theorem:
%
\begin{theorem}\label{thm:CLTL}
Assume $g \in L^2 (B(0, \delta))$ with
$\int_{B(0, \delta)} g(x)   \,dx \neq0$,
and let $N$ be a standard Gaussian random variable.
Then there is a constant $c_d$ only depending on the dimension $d$ such that
\[
d_W \biggl(\frac{F -\E F}{\sqrt{\V F}},N \biggr)\leq c_d
\lambda^{- 2} \max\bigl\{ 1, b(\delta)^{{3}/{2}} \bigr\}  \frac{  ( \int_{B(0, \d)} g( x)^4   \,dx  )^{1/2}} {
( \int_{B(0, \d)} g( x)   \,dx  )^2}.
\]
\end{theorem}

\begin{pf}
We compute the bound from Theorem \ref{thm:CLTlocal2}. Lemma \ref
{le:Leta} yields
\begin{eqnarray*}
\tilde V = \frac{\| f_1\|^2}{\lambda^3} &=& \int_{W} \biggl(\int
_{B(0, \d)} g( x) \1(y+x \in W ) \,dx \biggr)^2 \,dy
\\
&\geq& \int_{(1 /2) W} \biggl(\int_{B(0,\delta)}g(x)
\,dx \biggr)^2 \,dy = 2^{-d} \biggl(\int_{B(0, \d)}
g( x) \,dx \biggr)^2
\end{eqnarray*}
and
\begin{eqnarray*}
\bigl\| f^2\bigr\|_{\theta}^2 &=& \frac1{16} \int
_W \int_{B(0, \d)} g( x)^4
\1(y+x\in W) \,dx \,dy \leq\frac1{16} \int_{B(0, \d)} g(
x)^4 \,dx.
\end{eqnarray*}
\upqed\end{pf}

As an example we consider the particular case $g=1$, where $L(\eta)$
reduces to the number of edges of the graph.
Then the expectation is of order $\lambda^2 \delta^d$.
Lemma \ref{le:Leta} and Theorem \ref{thm:CLTL} tell us that the
variance is of order $\max \{\lambda^{3} \delta^{2d},\lambda
^2\delta^d \}$
and that
\[
d_W \biggl(\frac{L-\E L}{\sqrt{\V L}},N \biggr) \leq \tilde{c}_d
\lambda^{-2} \delta^{-3d/2 } \max\bigl\{ 1, \lambda^{{3}/{2}}
\delta^{{3d}/{2}} \bigr\}
\]
with a constant $\tilde{c}_d \in\R$ only depending on $d$.
The right-hand side tends to zero if $ \lambda^{{4}/{3}} \Lambda
_d(B(0, \delta)) \to\infty$ as $\lambda\to\infty$.
In the maybe\vspace*{1pt} most natural case when $\lambda\Lambda_d(B(0, \delta))$
stays constant we have an order
$\lambda^{- 1/2}$ of convergence to the Gaussian distribution. A
central limit theorem without rate of convergence is a special case of
Theorem 3.9 in \cite{PenroseRandomGeometricGraphs}.

Similar results to Lemma \ref{le:Leta} and Theorem \ref{thm:CLTL} can
be obtained if the intensity measure is of the form $d \mu(x)=\lambda
f(x)   \,d \Lambda_d(x)$ with $\lambda\in\R^+$ and a density function~$f(x)$.

\section*{Acknowledgments}
The authors are indebted to two anonymous referees for many helpful
remarks, and thank Christoph Th\"ale for valuable hints and helpful discussions.




\printaddresses

\end{document}